\documentclass[12pt]{amsart}
\usepackage{amsmath,amssymb,amsfonts}
\begin{document}
\theoremstyle{plain}
\newtheorem{thm}{Theorem}[subsection]
\newtheorem{lem}[thm]{Lemma}
\newtheorem{cor}[thm]{Corollary}
\newtheorem{prop}[thm]{Proposition}
\newtheorem{rem}[thm]{Remark}
\newtheorem{defn}[thm]{Definition}
\newtheorem{ex}[thm]{Example}

\numberwithin{equation}{subsection}
\newcommand{\mc}{\mathcal}
\newcommand{\mb}{\mathbb}
\newcommand{\surj}{\twoheadrightarrow}
\newcommand{\inj}{\hookrightarrow}
\newcommand{\zar}{{\rm zar}}
\newcommand{\an}{{\rm an}} 
\newcommand{\red}{{\rm red}}
\newcommand{\codim}{{\rm codim}}
\newcommand{\rank}{{\rm rank}}
\newcommand{\Ker}{{\rm Ker \ }}
\newcommand{\Pic}{{\rm Pic}}
\newcommand{\Div}{{\rm Div}}
\newcommand{\Hom}{{\rm Hom}}
\newcommand{\im}{{\rm im}}
\newcommand{\Spec}{{\rm Spec \,}}
\newcommand{\Sing}{{\rm Sing}}
\newcommand{\Char}{{\rm char}}
\newcommand{\Tr}{{\rm Tr}}
\newcommand{\Gal}{{\rm Gal}}
\newcommand{\Min}{{\rm Min \ }}
\newcommand{\Max}{{\rm Max \ }}
\newcommand{\sA}{{\mathcal A}}
\newcommand{\sB}{{\mathcal B}}
\newcommand{\sC}{{\mathcal C}}
\newcommand{\sD}{{\mathcal D}}
\newcommand{\sE}{{\mathcal E}}
\newcommand{\sF}{{\mathcal F}}
\newcommand{\sG}{{\mathcal G}}
\newcommand{\sH}{{\mathcal H}}
\newcommand{\sI}{{\mathcal I}}
\newcommand{\sJ}{{\mathcal J}}
\newcommand{\sK}{{\mathcal K}}
\newcommand{\sL}{{\mathcal L}}
\newcommand{\sM}{{\mathcal M}}
\newcommand{\sN}{{\mathcal N}}
\newcommand{\sO}{{\mathcal O}}
\newcommand{\sP}{{\mathcal P}}
\newcommand{\sQ}{{\mathcal Q}}
\newcommand{\sR}{{\mathcal R}}
\newcommand{\sS}{{\mathcal S}}
\newcommand{\sT}{{\mathcal T}}
\newcommand{\sU}{{\mathcal U}}
\newcommand{\sV}{{\mathcal V}}
\newcommand{\sW}{{\mathcal W}}
\newcommand{\sX}{{\mathcal X}}
\newcommand{\sY}{{\mathcal Y}}
\newcommand{\sZ}{{\mathcal Z}}
\newcommand{\A}{{\mathbb A}}
\newcommand{\B}{{\mathbb B}}
\newcommand{\C}{{\mathbb C}}
\newcommand{\D}{{\mathbb D}}
\newcommand{\E}{{\mathbb E}}
\newcommand{\F}{{\mathbb F}}
\newcommand{\G}{{\mathbb G}}
\renewcommand{\H}{{\mathbb H}}
\newcommand{\I}{{\mathbb I}}
\newcommand{\J}{{\mathbb J}}
\newcommand{\M}{{\mathbb M}}
\newcommand{\N}{{\mathbb N}}
\renewcommand{\P}{{\mathbb P}}
\newcommand{\Q}{{\mathbb Q}}
\newcommand{\R}{{\mathbb R}}
\newcommand{\T}{{\mathbb T}}
\newcommand{\U}{{\mathbb U}}
\newcommand{\V}{{\mathbb V}}
\newcommand{\W}{{\mathbb W}}
\newcommand{\X}{{\mathbb X}}
\newcommand{\Y}{{\mathbb Y}}
\newcommand{\Z}{{\mathbb Z}}


\catcode`\@=11
%
%
\def\opn#1#2{\def#1{\mathop{\kern0pt\fam0#2}\nolimits}} 
\def\bold#1{{\bf #1}}%
\def\underrightarrow{\mathpalette\underrightarrow@}
\def\underrightarrow@#1#2{\vtop{\ialign{$##$\cr
 \hfil#1#2\hfil\cr\noalign{\nointerlineskip}%
 #1{-}\mkern-6mu\cleaders\hbox{$#1\mkern-2mu{-}\mkern-2mu$}\hfill
 \mkern-6mu{\to}\cr}}}
\let\underarrow\underrightarrow
\def\underleftarrow{\mathpalette\underleftarrow@}
\def\underleftarrow@#1#2{\vtop{\ialign{$##$\cr
 \hfil#1#2\hfil\cr\noalign{\nointerlineskip}#1{\leftarrow}\mkern-6mu
 \cleaders\hbox{$#1\mkern-2mu{-}\mkern-2mu$}\hfill
 \mkern-6mu{-}\cr}}}
%
%

%
\def\:{\colon}
\let\oldtilde=\tilde
\def\tilde#1{\mathchoice{\widetilde{#1}}{\widetilde{#1}}%
{\indextil{#1}}{\oldtilde{#1}}}
\def\indextil#1{\lower2pt\hbox{$\textstyle{\oldtilde{\raise2pt%
\hbox{$\scriptstyle{#1}$}}}$}}
\def\pnt{{\raise1.1pt\hbox{$\textstyle.$}}}
%

%
\let\amp@rs@nd@\relax
\newdimen\ex@
\ex@.2326ex
\newdimen\bigaw@
\newdimen\minaw@
\minaw@16.08739\ex@
\newdimen\minCDaw@
\minCDaw@2.5pc
\newif\ifCD@
\def\minCDarrowwidth#1{\minCDaw@#1}
\newenvironment{CD}{\@CD}{\@endCD}
\def\@CD{\def\A##1A##2A{\llap{$\vcenter{\hbox
 {$\scriptstyle##1$}}$}\Big\uparrow\rlap{$\vcenter{\hbox{%
$\scriptstyle##2$}}$}&&}%
\def\V##1V##2V{\llap{$\vcenter{\hbox
 {$\scriptstyle##1$}}$}\Big\downarrow\rlap{$\vcenter{\hbox{%
$\scriptstyle##2$}}$}&&}%
\def\={&\hskip.5em\mathrel
 {\vbox{\hrule width\minCDaw@\vskip3\ex@\hrule width
 \minCDaw@}}\hskip.5em&}%
\def\verteq{\Big\Vert&&}%
\def\noarr{&&}%
\def\vspace##1{\noalign{\vskip##1\relax}}\relax\let\amp@rs@nd@&\iffalse}\fi
 \CD@true\vcenter\bgroup\relax\let\\=\cr\iffalse}\fi\tabskip\z@skip\baselineskip20\ex@
 \lineskip3\ex@\lineskiplimit3\ex@\halign\bgroup
 &\hfill$\m@th##$\hfill\cr}
\def\@endCD{\cr\egroup\egroup}
%
\def\>#1>#2>{\amp@rs@nd@\setbox\z@\hbox{$\scriptstyle
 \;{#1}\;\;$}\setbox\@ne\hbox{$\scriptstyle\;{#2}\;\;$}\setbox\tw@
 \hbox{$#2$}\ifCD@
 \global\bigaw@\minCDaw@\else\global\bigaw@\minaw@\fi
 \ifdim\wd\z@>\bigaw@\global\bigaw@\wd\z@\fi
 \ifdim\wd\@ne>\bigaw@\global\bigaw@\wd\@ne\fi
 \ifCD@\hskip.5em\fi
 \ifdim\wd\tw@>\z@
 \mathrel{\mathop{\hbox to\bigaw@{\rightarrowfill}}\limits^{#1}_{#2}}\else
 \mathrel{\mathop{\hbox to\bigaw@{\rightarrowfill}}\limits^{#1}}\fi
 \ifCD@\hskip.5em\fi\amp@rs@nd@}
\def\<#1<#2<{\amp@rs@nd@\setbox\z@\hbox{$\scriptstyle
 \;\;{#1}\;$}\setbox\@ne\hbox{$\scriptstyle\;\;{#2}\;$}\setbox\tw@
 \hbox{$#2$}\ifCD@
 \global\bigaw@\minCDaw@\else\global\bigaw@\minaw@\fi
 \ifdim\wd\z@>\bigaw@\global\bigaw@\wd\z@\fi
 \ifdim\wd\@ne>\bigaw@\global\bigaw@\wd\@ne\fi
 \ifCD@\hskip.5em\fi
 \ifdim\wd\tw@>\z@
 \mathrel{\mathop{\hbox to\bigaw@{\leftarrowfill}}\limits^{#1}_{#2}}\else
 \mathrel{\mathop{\hbox to\bigaw@{\leftarrowfill}}\limits^{#1}}\fi
 \ifCD@\hskip.5em\fi\amp@rs@nd@}
%
%
\newenvironment{CDS}{\@CDS}{\@endCDS}
\def\@CDS{\def\A##1A##2A{\llap{$\vcenter{\hbox
 {$\scriptstyle##1$}}$}\Big\uparrow\rlap{$\vcenter{\hbox{%
$\scriptstyle##2$}}$}&}%
\def\V##1V##2V{\llap{$\vcenter{\hbox
 {$\scriptstyle##1$}}$}\Big\downarrow\rlap{$\vcenter{\hbox{%
$\scriptstyle##2$}}$}&}%
\def\={&\hskip.5em\mathrel
 {\vbox{\hrule width\minCDaw@\vskip3\ex@\hrule width
 \minCDaw@}}\hskip.5em&}
\def\verteq{\Big\Vert&}
\def\novarr{&}
\def\noharr{&&}
\def\SE##1E##2E{\slantedarrow(0,18)(4,-3){##1}{##2}&}
\def\SW##1W##2W{\slantedarrow(24,18)(-4,-3){##1}{##2}&}
\def\NE##1E##2E{\slantedarrow(0,0)(4,3){##1}{##2}&}
\def\NW##1W##2W{\slantedarrow(24,0)(-4,3){##1}{##2}&}
\def\slantedarrow(##1)(##2)##3##4{%
\thinlines\unitlength1pt\lower 6.5pt\hbox{\begin{picture}(24,18)%
\put(##1){\vector(##2){24}}%
\put(0,8){$\scriptstyle##3$}%
\put(20,8){$\scriptstyle##4$}%
\end{picture}}}
\def\vspace##1{\noalign{\vskip##1\relax}}\relax\let\amp@rs@nd@&\iffalse}\fi
 \CD@true\vcenter\bgroup\relax\let\\=\cr\iffalse}\fi\tabskip\z@skip\baselineskip20\ex@
 \lineskip3\ex@\lineskiplimit3\ex@\halign\bgroup
 &\hfill$\m@th##$\hfill\cr}
\def\@endCDS{\cr\egroup\egroup}
%
\newdimen\TriCDarrw@
\newif\ifTriV@
\newenvironment{TriCDV}{\@TriCDV}{\@endTriCD}
\newenvironment{TriCDA}{\@TriCDA}{\@endTriCD}
\def\@TriCDV{\TriV@true\def\TriCDpos@{6}\@TriCD}
\def\@TriCDA{\TriV@false\def\TriCDpos@{10}\@TriCD}
\def\@TriCD#1#2#3#4#5#6{%
\setbox0\hbox{$\ifTriV@#6\else#1\fi$}
\TriCDarrw@=\wd0 \advance\TriCDarrw@ 24pt
\advance\TriCDarrw@ -1em
\def\SE##1E##2E{\slantedarrow(0,18)(2,-3){##1}{##2}&}
\def\SW##1W##2W{\slantedarrow(12,18)(-2,-3){##1}{##2}&}
\def\NE##1E##2E{\slantedarrow(0,0)(2,3){##1}{##2}&}
\def\NW##1W##2W{\slantedarrow(12,0)(-2,3){##1}{##2}&}
\def\slantedarrow(##1)(##2)##3##4{\thinlines\unitlength1pt
\lower 6.5pt\hbox{\begin{picture}(12,18)%
\put(##1){\vector(##2){12}}%
\put(-4,\TriCDpos@){$\scriptstyle##3$}%
\put(12,\TriCDpos@){$\scriptstyle##4$}%
\end{picture}}}
\def\={\mathrel {\vbox{\hrule
   width\TriCDarrw@\vskip3\ex@\hrule width
   \TriCDarrw@}}}
\def\>##1>>{\setbox\z@\hbox{$\scriptstyle
 \;{##1}\;\;$}\global\bigaw@\TriCDarrw@
 \ifdim\wd\z@>\bigaw@\global\bigaw@\wd\z@\fi
 \hskip.5em
 \mathrel{\mathop{\hbox to \TriCDarrw@
{\rightarrowfill}}\limits^{##1}}
 \hskip.5em}
\def\<##1<<{\setbox\z@\hbox{$\scriptstyle
 \;{##1}\;\;$}\global\bigaw@\TriCDarrw@
 \ifdim\wd\z@>\bigaw@\global\bigaw@\wd\z@\fi
 \mathrel{\mathop{\hbox to\bigaw@{\leftarrowfill}}\limits^{##1}}
 }
 \CD@true\vcenter\bgroup\relax\let\\=\cr\iffalse}\fi
 \tabskip\z@skip\baselineskip20\ex@
 \lineskip3\ex@\lineskiplimit3\ex@
 \ifTriV@
 \halign\bgroup
 &\hfill$\m@th##$\hfill\cr
#1&\multispan3\hfill$#2$\hfill&#3\\
&#4&#5\\
&&#6\cr\egroup%
\else
 \halign\bgroup
 &\hfill$\m@th##$\hfill\cr
&&#1\\%
&#2&#3\\
#4&\multispan3\hfill$#5$\hfill&#6\cr\egroup
\fi}
\def\@endTriCD{\egroup}

\title[An algebraic Riemann-Roch formula for flat bundles]{ A
remark on an algebraic Riemann-Roch formula for flat bundles }  
\author{ H\'el\`ene Esnault } 
\address{ Universit\"at Essen, FB6 Mathematik, 45 117 Essen, Germany}
\email{ esnault@uni-essen.de}
\thanks{ This work has been partly supported by the DFG Forschergruppe
''Arithmetik und Geometrie''}
\maketitle
\setcounter{section}{-1} 
\section{Introduction}
Let $X$ be a smooth variety defined over an algebraically closed
field $k$, and let $(E, \nabla)$ be a bundle with an integrable
connection. Then $(E, \nabla)$ carries algebraic classes $c_n
(E, \nabla)$ in the subgroup $\H^n (X, \Omega^{\infty} \sK_n)$
of the group of algebraic differential characters $AD^n (X)$
consisting of the classes mapping to 0 in $H^0 (X,
\Omega^{2n})$. These classes lift the Chern classes $c_n (E,
\nabla) \in CH^n (X)$, the algebraic Chern-Simons invariants
$w_n (E, \nabla) \in H^0 (X, \sH^{2n-1}_{DR})$ for $n \geq 2$,
as well as the analytic secondary invariants $c^{\an}_{n} (E, \nabla)
\in H^{2n-1} (X_{\an}, \C / \Z (n))$ when $k = \C$ (see
\cite{EII}, \cite{BE}, \cite{E} and section 1). 

Let now $f: X \to S$ be a smooth projective morphism. Then the
de Rham cohomology sheaves $R^j f_* (\Omega^{\bullet}_{X/S}
\otimes E, \nabla)$ carrie the (flat) Gau{\ss}-Manin connection,
still denoted by $\nabla$. Therefore one can ask for a
Riemann-Roch formula relating
$$
c_n (\sum (-1)^j (R^j f_* (\Omega^{\bullet}_{X/S} \otimes E ,
\nabla), \nabla))
$$
and $c_m (E, \nabla)$, as the classes $c_n (E, \nabla)$ verify the
Whitney product formula for exact sequences of bundles with
compatible (flat) connections. This formula should be compatible
with Riemann-Roch-Grothendieck formula in $CH^n (S)$, and with
Bismut-Lott and Bismut formula in $H^{2n-1} (X_{\an} , \C /\Z
(n))$ when $k = \C$ (\cite{Bi}, \cite{BiL}). In this note, we
propose an answer in the case $X = Y \times S$ and $f$ is the
projection: 

\begin{thm} \label{thm}
Let $Y$ be a smooth projective variety of dimension $d$ and $S$
be a smooth variety. Let $(E, \nabla)$ be a bundle with a flat
connection on $X = Y \times S$. Then 
$$
c_n (\sum_j (-1)^j R^j f_* (\Omega^{\bullet}_{X/S} \otimes E,
\nabla), \nabla) = (-1)^d f_* c_n (E, \nabla) |_{c_d
(\Omega^{1}_{X/S})} 
$$
\end{thm}
\noindent
where $f$ is the projection to $S$, where the right hand side
means that one takes a zero cycle $\Sigma = \sum_i m_i p_i \subset Y$
representing $c_d (\Omega^{1}_{Y})$, and 
$$
(-1)^d f_* c_n (E, \nabla) |_{c_d (\Omega^{1}_{X/S})} = \sum_{i}
m_i c_n (E, \nabla) |_{ \{ p_i\} \times S} \in \H^n (S,
\Omega^{\infty} \sK_n)
$$
(see section 1). 

For $n =1$, $\H^1 (S, \Omega^{\infty} \sK_1)$ is the group of
isomorphism class of rank one bundles with an integrable
connection. Thus the formula for the determinant bundle as a
flat bundle is similar to Deligne-Laumon formula for the
determinant bundle of a $\ell$-adic sheaf $E$ over $X$ with $S =
{\rm Spec} \ \bar{\F}_p$, $(p, \ell) = 1$ (\cite{D}, \cite{L}). But
here, we treat only the case where $f$ is split, which does not
occur in the $\ell$-adic content. However, we don't know what would be a
corresponding formula in higher codimension in the arithmetic
case. 

The method used is inspired by the work of Hitchin and Simpson.
We deform 
$$ 
\sum_j (-1)^j (R^j f_* (\Omega^{\bullet}_{X/S} \otimes E ,
\nabla ), \nabla)
$$ 
to the alternate sum of a flat structure
$$
(R^j f_* (\Omega^{\bullet}_{X/S} \otimes E, \alpha), \nabla_i)
$$
defined on the cohomology of a Higgs bundle, where the Higgs
structure is defined by a form $\alpha \in H^0 (Y,
\Omega^{1}_{Y})$, if $H^0 (Y, \Omega^{1}_{Y}) \neq 0$. If not,
one has to introduce poles (see section 2). 

\subsection{Acknowledgements:}

I started discussing on possible Riemann-Roch formulae with
Spencer Bloch. The motivation coming from \cite{D}, \cite{L} is
due to him. I thank him for his generosity. I thank G\'erard
Laumon and Domingo Toledo for explaining me parts of their work.
Finally I thank Eckart Viehweg for discussions on the use of
logarithmic forms and for his encouragements. 

\section{Trace of algebraic differential characters} 

\subsection{Algebraic differential characters} \label{subsec:AD}

Let $X$ be a smooth variety, $D \subset X$ be a normal crossing
divisor, $\sK_n$ be the Zarisky sheaf image of the Zarisky sheaf
of Milnor $K$-theory in $K^{M}_{n} (k(x))$. We recall \cite{E},
section 2, that the group 
$$
AD^n (X, D) : = \H^n (X, \sK_n \> d \log >> \Omega^{n}_{X} (\log
\ D) \to \ldots \to \Omega^{2n-1}_{X} (\log \ D))
$$
has a commutative product
$$
AD^m (X, D) \times AD^n (X, D) \to AD^{m+n} (X,D),
$$
respecting the subgroup 
\begin{multline*}
\H^n (X, \Omega^{\infty}_{X} (\log \ D) \sK_n ): = \\
\H^n (X, \sK_n \> d \log >> \Omega^{n}_{X} (\log \ D) \to \ldots ) \\
= {\rm Ker} \ (AD^n (X, D) \to H^0 (X, \Omega^{2n}_{X} (\log \ D))
\end{multline*}
and compatible with the products in $\sK_n$ and $\Omega^{\geq n}_{X} (\log 
\ D)$. A bundle $(E, \nabla)$ with a $\Omega^{1}_{X} (\log \ D)$ connection
has functorial and additive classes $c_n (E, \nabla) \in AD^n
(X, D)$, lying in $\H^n (X, \Omega^{\infty}_{X} (\log \ D)
\sK_n)$ when $\nabla^2 = 0$. 

\begin{lem} \label{L} 

Let $X = \P^1 \times S$, with $S$ smooth, $B \subset \P^1$ be a
divisor, $D = B \times S$, $p_i , i = 1,2$ be the projections of
$X$ to $\P^1$ and $S$. Then one has a direct sum decomposition
\begin{multline*}
p^{*}_{2} \oplus p^{*}_{2} \otimes p^{*}_{1} : \\
\H^n (S, \Omega^{\infty}_{S} \sK_n) \oplus \H^{n-1} (S,
\Omega^{\infty}_{S} \sK_{n-1}) \otimes \H^1 (\P^1 ,
\Omega^{\infty}_{\P^1} (\log \ B) \sK_1) \\
\to \H^n (X, \Omega^{\infty}_{X} (\log \ D) \sK_n). 
\end{multline*}
\end{lem}

\begin{proof} One has the K\"{u}nneth formulae
\begin{multline*}
H^{\ell} (X, \sK_n) = \\
p^{*}_{2} H^{\ell} (S, \sK_n) \oplus p^{*}_{2} H^{\ell-1} (S, \sK_{n-1})
\cup p^{*}_{1} H^1 (\P^1, \sK_1) 
\end{multline*}
\begin{multline*}
\H^{\ell} (X, \Omega^{\geq n}_{X} (\log \ D)) = \\p^{*}_{2} \H^{\ell} (S,
\Omega^{\geq n}_{S} ) 
\oplus p^{*}_{2} \H^{\ell-1} (S, \Omega^{\geq n -1}_{S} ) \cup
p^{*}_{1} \H^0 (\P^{1} , \Omega^{1}_{\P^1} (\log \ B)) \\
\mbox{if} \ \ B \neq \phi 
\end{multline*}
\begin{multline*}
\H^{\ell} (X, \Omega^{\geq n}_{X} (\log \ D)) = \\
= p^{*}_{2} \H^{\ell} (S, \Omega^{\geq n}_{S}) 
\oplus p^{*}_{2} \H^{\ell-2} (S, \Omega^{\geq n-1}_{S}) \cup
p^{*}_{1} H^0 (\P^1 , \Omega^{1}_{\P^1}) \\
\mbox{if} \ \ B = \phi. 
\end{multline*}
Moreover, in the long exact cohomology sequence associated to
the short exact sequence 
$$
0 \to \Omega^{\geq n}_{X} (\log \ D) [n-1] \to
\Omega^{\infty}_{X} (\log \ D) \sK_n \to \sK_n \to 0 
$$
the map
$$
H^{\ell} (X, \sK_n) \> d \ \text{log} >> \H^{\ell+n} (X, \Omega^{\geq n}_{X}
(\log \ D))
$$
respects this direct sum composition. (If $B \neq \phi$, the
term 
$$
p^{*}_{2} H^{\ell-1} (S, \sK_{n-1} ) \cup p^{*}_{1} H^1 (\P^1 ,
\sK_1) 
$$
maps to zero). This shows the lemma. 
\end{proof}

\subsection{Trace} \label{subsec:tr} 
\begin{prop} \label{P} 

Let $X = Y \times S$, with $Y$ smooth projective, $S$ smooth, and
$f: X \to S$ be the projection. Let $\Sigma = \sum m_i p_i$ be a
zero cycle in $Y$, and $(E, \nabla)$ be a bundle with an
integrable connection on $X$. Then 
$$ 
f_* c_n (E, \nabla) |_{\Sigma \times S} : = \sum m_i c_n ((E,
\nabla) |_{ \{ p_i \} \times S} )
\in \H^n (S, \Omega^{\infty}_{S} \sK_n)
$$
does not depend on the choice of the representative $\Sigma$ in
its equivalence class $[\Sigma ] \in CH_0 (S)$. 
\end{prop}

\begin{proof} Let $\Sigma' = \sum m_i p'_i$ be another choice. Then
there are rational functions $f_i $ on curves $C_i \subset Y$
such that $\Sigma - \Sigma' = \sum$ div $f_i$. Therefore it is
sufficient to prove the following: let $\nu: C \to Y$ be the
normalization of an irreducible curve $\nu (C) \subset X$, and
$\varphi: C \to \P^1$ be a non trivial rational function on
$C$. Then 
$$
c_n ((\nu \times \text{id}_S )^* (E, \nabla) |_{\varphi^{-1} (0) \times
S}) = c_n ((\nu \times \text{id}_S)^* (E, \nabla) |_{\varphi^{-1}
(\infty) \times S}). 
$$
Let $(\sE , D) : = (\nu \times \text{id}_X)^* (E, \nabla)$. Let $B
\subset \P^1$ be the ramification locus of $\varphi$, $\pi =
\varphi \times \text{id}_S : C \times S \to \P^1 \times S$. Then $\pi_*
(\sE , D)$ is a bundle with an integrable connection
with logarithmic poles along
$D = B \times S$, and has classes 
$$ 
c_n ( \pi_* (\sE, D)) \in \H^n (\P^{1} \times S,
\Omega^{\infty}_{X} (\log \ D) \sK_n).
$$
For $t \notin B$, then 
$$
c_n (\pi_* (\sE, D) |_{ \{ t\} \times S}) = \sum_{s \in
\varphi^{-1}(t)} c_n ( (E, D) |_{ \{ s \} \times S})= c_n (\sE,
D) |_{\varphi^{-1} (t) \times S} .
$$
For $t \in B$, the residue map 
$$
{\rm res}_{\{ t \} \times S} : p^{*}_{2} \Omega^{n-1}_{S} 
\otimes p^{*}_{1} \Omega^{1}_{\P^1} (\log \ B) \to
\Omega^{n-1}_{\{ t \} \times S} 
$$
verifies 
$$
{\rm res}_{\{ t\} \times S} \otimes _{\sO_{\P^1
\times S}} \sO_{\{ t \} \times S} = {\rm Identity}_{\{ t\} \times
S} .
$$
Therefore there is a canonical splitting 
$$
\Omega^{n}_{X} (\log \ D)|_{\{ t \} \times S} = \Omega^{n}_{S}
\oplus \Omega^{n-1}_{\{ t \} \times S} 
$$
Thus the connection 
$$
\pi_* D : \pi_* \sE \to \Omega^{1}_{\P^1 \times S} (\log \ D)
\otimes \pi_* \sE
$$
restricts to 
$$
\pi_* D |_{\{ t \} \times S} : \pi_* \sE |_{\{ t \} \times S} \to
(\Omega^{1}_{S} \oplus \sO_S) \otimes \pi_* \sE |_{\{ t \}
\times S} 
$$
defining by projection the connection 
$$
\overline{\pi_* D} |_{\{ t \} \times S} : \pi_* \sE |_{ \{ t \}
\times S} \to \Omega^{1}_{S} \otimes (\pi_* \sE) |_{\{ t \}
\times S}.
$$
The integrability of $\pi_* D |_{\{ t\} \times S} $ implies the
integrability of the genuine connection $\overline{\pi_* D}
|_{\{ t \} \times S}$, and 
$$
c_n (( \sE, D) |{\varphi^{-1} (t) \times S}) = c_n ( \pi_* \sE
|_{ \{ t \} \times S} , \overline{\pi_* D} |_{ \{ t\} \times
S}). 
$$
Now by \cite{L}, 
$$
c_n (\pi_* (\sE, D)) = p^{*}_{2} a + \sum_{i} p^{*}_{2} b_i \cup
p^{*}_{1} c_i
$$
where $a \in \H^n (S, \Omega^{\infty}_{S} \sK_n)$, and the preceding
discussion shows that 
$$
c_n ((\sE , D) |_{\varphi^{-1} (t) \times S} ) = a \in \H^n (S,
\Omega^{\infty}_{S} \sK_n)
$$
whether $t \in B$ or $t \notin B$.
\end{proof} 

\section{Proof of the theorem} 

\subsection{Notations} 

Let $H_0 , \ldots , H_N$ be 
restriction of the coordinate hyperplanes $H'_0, \ldots , H'_N$ 
in an embedding $ Y \subset \P^N$. Since $\Omega^1_{\P^N}(\log
H'_0 \cup \ldots H'_N) \cong \oplus_1^N \sO_{\P^N}$, the sheaf 
$\Omega^1_Y(\log H_0 \cup \ldots H_N)$ is globally generated.
We denote by
$H^{(\ell)}$ the normalization of the $\ell$ by $\ell$ intersections of
the $H_{j}$, $H^{(0)} = X$, $H^{(\delta)} = \emptyset$ for $\delta >
d$, by $H$ the union of the $H_j$, by 
$$
\Omega^{a}_{H^{(\ell)}} (\log \ H^{(\ell+1)})
$$
the sheaf of a form on $H^{(\ell)}$ with logarithmic poles along
the $(\ell+1)$ by $(\ell+1)$ intersections. One has the following
resolution of the de Rham complex

\begin{equation} \label{E}
\Omega^{\bullet}_{X/S} \to \Omega^{\bullet}_{X/S} (\log \ (H
\times S)) \to \Omega^{\bullet -1}_{H^{(1)} \times S/S} (\log
(H^{(2)} \times S)) \to \ldots 
\end{equation}
$$
\to \Omega^{\bullet -\ell}_{H^{(\ell)} \times S/S} (\log \ (H^{(\ell+1)}
\times S)) \to \ldots 
$$
$$
\to \Omega^{\bullet -d}_{H^{(d)} \times S/S} \to 0.
$$ 

\subsection{Proof} 

The resolution \ref{E} is compatible with the Gau{\ss}-Manin
connection, as
$$
H^{(\ell)} \times S
$$
is dominant over $S$. 

Therefore in the $K$ group $K(S,$ flat) of bundles on $S$
with an integrable connection, one has
\begin{multline*}
\sum_j (-1)^j R^j f_* (\Omega^{\bullet}_{X/S} \otimes E,
\nabla)\\ 
= \sum (-1)^j R^j f_* (\Omega^{\bullet}_{X/S} (\log \ H \times
S) \otimes E, \nabla) \\
+ \sum (-1)^j R^j f_* (\Omega^{\bullet}_{H^{(1)} \times S/S}
(\log \ (H^{(2)} \times S)) \otimes E, \nabla) \\
+ \ldots + \sum (-1)^j R^j f_* (\Omega^{\bullet}_{H^{(\ell)} \times
S/S} (\log \ (H^{(\ell)} \times S )) \otimes E, \nabla) \\
+ \ldots + f_* (\Omega^{\bullet}_{H^{(d)} \times S/S} \otimes E,
\nabla). 
\end{multline*}
On the other hand 
\begin{multline*}
c_d (\Omega^{\bullet}_{X/S} ) = c_d (\Omega^{1}_{X/S} (\log \
H)) \\
- c_{d-1} (\Omega^{1}_{H^{\ell} \times S/S} (\log \ (H^{(2)} \times S)) \\
+ \ldots + (-1)^{\ell} c_{d-\ell} (\Omega^{1}_{H^{(\ell)} \times S/S} (\log
\ (H^{(\ell+1)} \times S )) \\
+ \ldots + (-1)^d [H^{(d)} \times S]
\end{multline*}
where $[H^{(d)} \times S]$ means the codimension $d$ cycle,
image of $H^{(d)} \times S$ in $Y \times S$. 

Therefore, one just has to prove the following formula
\begin{multline*} 
c_n ( \sum (-1)^j (R^j f_* (\Omega^{\bullet}_{X/S} (\log \ D
\times S) \otimes E) , \nabla), \nabla)) \\
= (-1)^d f_* c_n (E, \nabla) |_{c_d (\Omega^{1}_{X/S} (\log \ (H
\times S))} ).
\end{multline*} 
We denote by $\tau : \Omega^{1}_{X} (\log \ (H \times S)) \to
f^* \Omega^{1}_{S}$ the splitting of the one forms, which induces
a splitting 
$$
\tau : \Omega^{\bullet}_{X} (\log \ (H \times S) ) \to f^*
\Omega^{\bullet}_{S} 
$$
of the de Rham complex, where the differential $f^*
\Omega^{i}_{S} \to f^* \Omega^{i+1}_{S} $ is defined by $\tau d
\iota$, $\iota : f^* \Omega^{i}_{S} \to \Omega^{i}_{X} (\log \
(H \times S))$ being the natural embedding. This defines a $f^*
\Omega^{1}_{S}$ valued connection $\nabla_{\tau}$ on
$\Omega^{i}_{X/S} (\log \ (H \times S)) \otimes E$ by embedding
$\Omega^{i}_{X/S} (\log \ (H \times S)) \otimes E$ into
$\Omega^{i}_{X} (\log \ (H \times S)) \otimes E$ via the
splitting, then taking $\nabla$, then projecting onto the factor
$$
f^* \Omega^1_S \otimes \Omega^{1}_{X/S} (\log \ (H \times S)) \otimes E
$$
with the sign $(-1)^i$. The integrability condition $\nabla^2
=0$ then implies that
$$
\nabla_{\tau} \circ \nabla_{X/S} = \nabla_{X/S} \circ
\nabla_{\tau}, 
$$
where 
$$
\nabla_{X/S} : \Omega^{i}_{X/S} (\log \ (H \times S)) \otimes E
\to \Omega^{i+1}_{X/S} (\log \ (H \times S)) \otimes E.
$$
Taking cohomology defines a connection, still denoted by $\nabla_{\tau}$
$$
\nabla_{\tau} : R^j f_* \Omega^{i}_{X/S} (\log \ (H \times S))
\otimes E \to \Omega^{1}_{S} \otimes R^j f_* \Omega^{i}_{X/S}
(\log \ (H \times S) \otimes E.
$$
The integrability of $\nabla$ implies the integrability of
$\nabla_{\tau}$. Therefore in $K(S$, flat) one has 
\begin{multline*}
(R^j (f_* (\Omega^{\bullet}_{X/S} (\log \ (H \times S)) \otimes
E, \nabla) , \nabla) \\
= \oplus_i  R^{j-i} (f_* \Omega^{i}_{X/S} (\log \ (H \times S)
) \otimes E, \nabla_{\tau} ) , \nabla_{\tau}) 
\end{multline*}
and
\begin{multline*}
\sum_j (-1)^j (R^j f_* \Omega^{\bullet}_{X/S} ( \log \ (H \times
S) \otimes E, \nabla) , \nabla) \\
= \sum_{i,j} (-1)^{i+j} (R^j f_* \Omega^{i}_{X/S} (\log \ (H
\times S)) \otimes E, \nabla_{\tau}). 
\end{multline*}
Let $\alpha
\in H^0 (Y, \Omega^{1}_{Y} (\log \
H))$ be a non-trivial generic section. 
We still denote by $\alpha$ the corresponding form
$p^{*}_{1} \alpha \in H^0 (X, \Omega^{1}_{X/S} (\log \ (H \times
S)))$. Then, it defines a morphism 
$$
\alpha_{X/S} : \Omega^{i}_{X/S} (\log \ (H \times S)) \otimes E
\to \Omega^{i+1}_{X/S} (\log \ (H \times S) ) \otimes E
$$
by $\alpha_{X/S} (w \otimes e) = \alpha  \wedge w \otimes e$. As
$d \alpha =0$, one has 
$$
\alpha_{X/S} \circ \nabla_{\tau} = \nabla_{\tau} \circ
\alpha_{X/S}.
$$
Thus in $K$ ($S$, flat) one has 
\begin{multline*}
\sum_{i,j} (-1)^{i+j} (R^j f_* (\Omega^{i}_{X/S} (\log \ (H
\times S)) \otimes E, \nabla_{\tau} ), \nabla_{\tau}) \\
= \sum_j (-1)^j (R^j f_* (\Omega^{\bullet}_{X/S} (\log \ (H
\times S)) \otimes E, \alpha_{X/S} ), \nabla_{\tau}).
\end{multline*}
On the other hand, the complex
$$
(\Omega^{\bullet}_{X/S} (\log \ (H \times S)), \alpha_{X/S})
$$ 
is quasi-isomorphic to $\sO_{\Sigma \times S} [-d]$, where $\Sigma$ is
the zero set of $\alpha$, and one furthermore has a commutative
diagram of complexes 
$$
\begin{CD}
(\Omega^{\bullet}_{X/S} (\log \ (H \times S)) \otimes E, \alpha
) \> \nabla_{\tau} >> f^* \Omega^{1}_{S} \otimes
(\Omega^{\bullet}_{X/S} (\log \ (H \times S)) \otimes E, \alpha
) \\
\V VV \V VV \\
E|_{\Sigma \times S} [-d] \> \nabla |_{\Sigma \times S} >>
\Omega^{1}_{\Sigma \times S} \otimes E [-d]
\end{CD}
$$
This shows the equality in $K$ ($S$, flat)
\begin{multline*}
\sum_j (-1)^j (R^j f_* (\Omega^{\bullet}_{X/S} ( \log \ (H
\times S) \otimes E, \alpha) , \nabla_{\tau}) \\
= (-1)^d \oplus_{\sigma \in \Sigma} (E, \nabla) |_{\sigma \times
S} 
\end{multline*}
and finishes the proof. 

\bibliographystyle{plain}
\renewcommand\refname{References}
 
\end{document}